\documentclass[12pt]{article}
\newcommand{\eqdef}{\, =\kern -12.7pt\raise 6pt\hbox{{\tiny\textrm{def}}}\,\,}
\newtheorem{theorem}{Theorem}
\title{
There's plenty of time for evolution}
\author{ Herbert S. Wilf\\Department
of Mathematics, University of Pennsylvania\\Philadelphia, PA
19104-6395\\
{\small\tt <wilf@math.upenn.edu>} \and Warren J. Ewens\\
Department of Biology, University of Pennsylvania\\
Philadelphia, PA 19104-6018\\
{\small\tt <wewens@sas.upenn.edu>} }
\begin{document}
\maketitle
\begin{abstract}
Objections to Darwinian evolution are often based on the time required to carry out the necessary mutations. Seemingly, exponential numbers of mutations are needed. We show that such estimates ignore the effects of natural selection, and that the numbers of necessary mutations are thereby reduced to about $K\log{L}$, rather than $K^L$, where $L$ is the length of the genomic ``word,'' and $K$ is the number of possible ``letters'' that can occupy any position in the word. The required theory makes contact with the theory of radix-exchange sorting in theoretical computer science, and the asymptotic analysis of certain sums that occur there.
\end{abstract}
\section{Introduction}
The 2009 ``Year of Darwin'' has seen many
 welcome tributes to this great scientist, and
 re-affirmations of the validity of his theory of evolution
 by natural selection, though this validity is not universally accepted. One of the main objections that have been raised holds that there has not been enough time for all of the species complexity that we see to have evolved by random mutations. Our purpose here is to analyze this process, and our conclusion is that when one takes account of the role of natural selection in a reasonable way, there has been ample time for the evolution that we observe to have take place.

 \section{The calculations}

 Biological evolution is such a complex process
 that any attempt to describe it precisely in a way similar to the description
 of the dynamic processes in physics by mathematical methods is impossible. This does not mean,
 however, that arbitrary models of biological evolution are allowed. Any allowable model has
 to reflect the main features of evolution. Our main aims, discussed
 below, are to indicate why an evolutionary model often used to
 ``discredit'' Darwin, leading to the ``not enough time'' claim,  is inappropriate, and to find the mathematical properties
 of a more appropriate model.

 Before doing this we take up some other points. Evolution as a Darwinian-Mendelian
 process takes place via a succession of gene replacement processes, whereby a new
 ``superior'' gene arises by mutation  in the population and, by natural selection, steadily replaces
 the current gene. (We use here the word ``gene'' rather than the more technically accurate ``allele''.)
 It has recently been estimated \cite{yx} that a newborn human carries
 some de novo 100-200 base mutations. Only about five of these
 can be expected, on average, to arise in parts of the genome coding for
 genes or in regulatory regions. In a population admitting a million births in any year, we may expect
 something on the order of five million such de novo mutations, or about 250 per gene
 in a genome containing 20,000 genes.
There is then little problem about a supply of new mutations in any gene. However
only a small proportion of these can be expected to be favorable. We formalize this
in the calculations below.

We now
turn to the inappropriate evolutionary model referred
to above concerning the fixation of these genes in the population.
The incorrect argument runs along the following lines. Consider the replacement processes
needed  in order to
change each of the resident genes at $L$ loci  in a more primitive genome into those of a more
favorable, or advanced, gene. Suppose that at each such gene locus,  the argument runs,
the proportion of gene types (alleles) at that gene locus that are more favored than  the primitive type
is $K^{-1}$.   The probability that at all $L$ loci a more favored gene type  is obtained
in one round of evolutionary ``trials''
is  $K^{-L}$, a vanishingly small amount. When trials are carried out
sequentially over time,  an
exponentially large number of trials (of order $K^L$) would be needed in order to carry out
the complete transformation, and from this some have concluded that the evolution-by-mutation paradigm doesn't work because of lack of time.

But this argument in effect
assumes an ``in series'' rather than a more correct ``in parallel''  evolutionary
process. If a superior gene for (say) eye function has become fixed in a
population, it is not thrown out when a superior gene for (say) liver
function becomes fixed. Evolution is an ``in parallel'' process, with beneficial
mutations at one gene locus being retained after they become fixed in a population while  beneficial
mutations at other loci become fixed. In fact this statement is essentially the principle of natural selection.

The paradigm used in the incorrect argument is often formalized as follows.
Suppose that we are trying to find a specific unknown word of $L$ letters, each of the letters having
been chosen from an alphabet of $K$ letters. We want to find the word by means of
a sequence of rounds of guessing letters. A single round consists in guessing all
of the letters of the word by choosing, for each letter, a randomly chosen letter from
the alphabet. If the correct word is not found, a new sequence is guessed, and the procedure is continued
until the correct sequence is found. Under this paradigm the mean number
of rounds of guessing until the correct sequence is found is indeed $K^{L}$.

But a more appropriate model is the following. After guessing each of the letters,
we are told which (if any) of the guessed
letters are correct, and then those letters are retained.
The second round of guessing is applied only for the incorrect letters that remain after this
first round, and so forth. This procedure mimics the ``in parallel''
evolutionary process. The question concerns the statistics of the number of
rounds needed to guess all of the letters of the word successfully.
Our main result is
\begin{theorem}
The mean number of rounds that are necessary to guess all of the letters
 of an $L$ letter word, the letters coming from an alphabet of $K$ letters, is
\begin{equation}
= \frac{\log{L}}{\log{(\frac{K}{K-1})}}+\beta(L)+O(L^{-1})\qquad\quad(L\to\infty).
\end{equation}
with $\beta(L)$ being the periodic function of $\log{L}$ that is given by
 eq. (\ref{eq:beta3}) below. The function $\beta(L)$ oscillates
  within a range
  which for $K\ge 2$, is never larger than $.000002$ about the first two terms
  on the right-hand side of equation (\ref{eq:beta3}).
\end{theorem}

For example, if we are using a $K=40$ letter alphabet, and a word of
length 20,000 letters, then the number of possible words is about $10^{34,040}$, but
our theorem shows that a mean of only about
\[\frac{\log{20,000}}{\log{(\frac{40}{39})}}\approx 390\]
rounds of guessing will be needed, where each round consists of one pass through
the entire as-yet-unguessed word.

The central feature of this result lies in the logarithmic terms in the above expression. Even if $L$
is very large, log $L$ is (for values of $L$ arising in practice in any genome) in practice
manageable. The inappropriate arguments referred to above lead to the value $10^{34,040}$,
and arise because of the incorrect ``in series'' rather than the correct ``in parallel''
implicit assumption about the nature of genetic evolution.

We have chosen the numerical values in this example to reflect the biological evolutionary process.
The value 20,000 represents the number of gene loci  in the genome at which replacement
processes are to take place. The value $K = 40$ is arrived at by using the value 250
found above for the number of de novo mutations per gene locus per year and a rough estimate
that only one mutation in 10,000 is selectively favored over the resident gene type.  In
practice further modifications are needed to the calculations since, because of stochastic
events, only a proportion of selectively favored new mutations become fixed in a population.

However, the force of our result does not depend on the numerical values that one assigns to $K$ and $L$. The fact is that with the parallel model, i.e., taking account of natural selection, the number of rounds of mutations that are needed to change the complete genome to its desirable form are only about $K\log{L}$, instead of the hugely exponential $K^L$ which would result from the serial model.

\section{The analysis}
The probability that the first letter of the word will be correctly guessed  in at
most $r$ rounds of guessing is
\[1-\left(1-\frac{1}{K}\right)^r,\]
so the probability that all $L$ letters of the given word will be guessed correctly
 in $\le r$ rounds is
\[\left(1-\left(1-\frac{1}{K}\right)^r\right)^L.\]
Thus the mean number of rounds that will be needed to guess all of the letters
 of the word is
\[\sum_{r=1}^{\infty}r\left\{\left(1-\left(1-\frac{1}{K}\right)^r\right)^L
-\left(1-\left(1-\frac{1}{K}\right)^{r-1}\right)^L\right\},\]
which is simply the mean of the maximum of $L$ independently and identically distributed
(iid)  geometric random variables.

This infinite sum can be transformed into a finite sum because
\begin{eqnarray*}
\sum_rr\left\{(1-x^r)^L-(1-x^{r-1})^L\right\}&=&\sum_rr\sum_j\left\{{L\choose j}(-1)^jx^{rj}-{L\choose j}(-1)^jx^{(r-1)j}\right\}\\
&=&\sum_{j=0}^L{L\choose j}\frac{(-1)^{j+1}}{1-x^j}.
\end{eqnarray*}

Consequently the mean number of rounds needed to guess all $L$ letters is
\begin{equation}\label{eq:avg}\alpha(L)\eqdef\sum_{j=1}^L{L\choose j}\frac{(-1)^{j+1}}{1-(1-1/K)^j}=1+\sum_{j=1}^L(-1)^{j+1}{L\choose j}\frac{1}{X^j-1}.\qquad(X=K/(K-1))\end{equation}
The appearance of this latter sum in the current context is somewhat surprising. It is one which is well known to theoretical computer scientists, and it arises there in connection with radix-exchange sorting.

To find the asymptotic behavior of this sum, we note that the behavior of the
following sum is known, and can be found in Exercise 50 of section 5.2.2 of \cite{kn}:
\begin{equation}\label{eq:usum}
U_{m,n}=\sum_{k\ge 2}{n\choose k}\frac{(-1)^k}{m^{k-1}-1}, \qquad(m>1).
\end{equation}
The result in \cite{kn}, due to N.G. deBruijn, is that
\begin{equation}\label{eq:usum2}
U_{m,n}=n\log_m{n}+n\left(\frac{\gamma-1}{\log{m}}-\frac12+f_{-1}(n)\right)
+\frac{m}{m-1}-\frac{1}{2\log{m}}-\frac12f_{1}(n)+O(n^{-1}),
\end{equation}
where $\gamma=.57721..$ is Euler's constant and
\[f_{s}(n)=\frac{2}{\log{m}}\sum_{k\ge 1}\Re(\Gamma(s-2\pi ik/\log{m})\exp{(2\pi ik\log_m(n))}).\]
These $f_{s}(n)$'s are bounded functions, and in fact, they are evidently periodic
of period 1 in $\log_m(n)$.

To relate our sum to Knuth's, we have
\[\alpha(L)=1+U_{X,L+1}-U_{X,L}.\qquad(X=K/(K-1)).\]
We note in passing that this shows that the mean of $L$ iid geometric random variables
 is, quite generally, simply related to the quantities $U_{m,n}$, which had
 previously been encountered in connection with radix-exchange sorting, and whose
 (notoriously difficult) asymptotic behavior had been found as a result of that
 connection.

After doing the subtraction we obtain
\begin{eqnarray}
\alpha(L)&=&\frac{\log{L}}{\log{X}}+L(f_{-1}(L+1)-f_{-1}(L))+\frac12
+\frac{\gamma}{\log{X}}+\frac12(f_{1}(L+1)-f_{1}(L))+O(L^{-1})\nonumber\\
&\eqdef&\frac{\log{L}}{\log{X}}+\beta(L)+O(L^{-1}),\label{eq:beta2}
\end{eqnarray}
where again $X=K/(K-1)$, and we have written
\begin{equation}\label{eq:beta}
\beta(L)= L(f_{-1}(L+1)-f_{-1}(L))+\frac12
+\frac{\gamma}{\log{X}}+\frac12(f_{1}(L+1)-f_{1}(L)).
\end{equation}
But we have
\[L(f_{-1}(L+1)-f_{-1}(L))=\frac{2}{\log{m}}\sum_{k\ge 1}\Re\left(\Gamma(-1-2\pi ik/\log{m})\exp{(2\pi ik\log_m{L})}\frac{2\pi ik}{\log{m}}\right)+O(L^{-1}),\]
whereas
\[f_{1}(L+1)-f_{1}(L)=O(L^{-1}).\]
Therefore for $\beta(L)$ in (\ref{eq:beta2}) we can take
\begin{equation}\label{eq:beta3}
\beta(L)=\frac12
+\frac{\gamma}{\log{X}}+\frac{2}{\log{X}}\sum_{k\ge 1}\Re\left(\Gamma(-1-2\pi ik/\log{X})\exp{(2\pi ik\log_X{L})}\frac{2\pi ik}{\log{X}}\right).
\end{equation}

Therefore the mean number of rounds that are necessary to guess all of the letters
 of an $L$ letter word, the letters coming from an alphabet of $K$ letters, is given
  exactly by (\ref{eq:avg}), and is asymptotically
\begin{equation}\label{eq:beta4}
= \frac{\log{L}}{\log{(\frac{K}{K-1})}}+\beta(L)+O(L^{-1})\qquad\quad(L\to\infty).
\end{equation}
with $\beta(L)$ being the periodic function of period 1 in $\log_X{L}$ that is given by
 (\ref{eq:beta3}) and $X=K/(K-1)$. If terms of order $L^{-1}$ and the small constant 0.000002
 mentioned in the statement of the theorem are ignored, this is
 \begin{equation}\label{eq:beta455}
\frac{\log{L}+\gamma }{\log{(\frac{K}{K-1})}}+\frac12.
\end{equation}

 We conclude with a comment on the oscillatory behavior of this mean, as revealed
 by the exact expression in (\ref{eq:beta3}).  In probability theory the asymptotic behavior of a maximum of several iid random variables is often found by ``sandwiching'' the discrete random variable between two continuous
 random variables whose asymptotic behavior is known. In the case of geometric random variables the
 appropriate continuous random variables to be used for this purpose have negative exponential
 distributions. This sandwiching procedure has been used frequently and leads to the expression (\ref{eq:beta455}),
 but does not lead to the more precise oscillatory behavior exhibited in the expression (\ref{eq:beta3}).

\end{document}